# Growth of the Higgs field for solutions to the Kapustin-Witten equations on $\mathbb{R}^4$

Clifford Henry Taubes†

Department of Mathematics
Harvard University
Cambridge, MA 02138

chtaubes@math.harvard.edu

ABSTRACT: The Kapustin-Witten equations on $\mathbb{R}^4$ are equations for a pair of connection on the product principle SU(2) bundle and 1-form with values in the product Lie algebra bundle. The 1-form is the Higgs field. A dichotomy is proved to the effect that either the averaged norm of the Higgs field on large radius spheres grows faster than a power of the radius, or its 1-form components everywhere pairwise commute.

†Supported in part by the National Science Foundation

## 1. Introduction

To set the stage for what is to come, let SU(2) denote the Lie group of $2 \times 2$ unitary matrices with determinant equal to 1. Let $\mathfrak{su}(2)$ denote its Lie algebra, and let $\mathbb{V}$ denote an auxilliary, finite dimensional vector space over $\mathbb{R}$ with an inner product. Supposing that $n \geq 2$, let $(A, a)$ denote a pair of connection on the product SU(2) bundle over $\mathbb{R}^n$ and section of the associated vector bundle with fiber $\mathbb{V} \otimes \mathfrak{su}(2)$. The pair $(A, a)$ is constrained to obey

$$\nabla_A^\dagger \nabla_A a + [a_c, [a, a_c]] = 0$$

(1.1)

where $\nabla_A$ denotes the covariant derivative defined by A, where $\nabla_A^\dagger$ denotes the formal $L^2$ adjoint of this covariant derivative, and where $\{a_c\}_{c=1,\ldots,\dim(\mathbb{V})}$ are the components of $a$ with respect to any given orthonormal frame for the vector space $\mathbb{V}$. The brackets in (1.1) denote commutators; and repeated indices are summed. (This repeated index convention holds in subsequent equations also.)

Let $\kappa$ denote the non-negative function on $[0, \infty)$ whose square is given by

$$\kappa^2(r) = \frac{1}{r^{n-1}} \int_{\partial B_r} |a|^2$$

(1.2)

with $B_r$ denoting the radius r ball about the origin and with $\partial B_r$ denoting its boundary sphere. If $a$ is not identically zero, then the function $\kappa$ is non-decreasing on $(0, \infty)$ (see the upcoming Equation (2.1).) As a consequence, there are no interesting solutions where $\kappa$ has limit zero as $r \to \infty$. The theorem that follows says more about how the large r asymptotics of $\kappa$ constrain the solutions to (1.1). This theorem asserts in effect a dichotomy concerning the growth of the function $\kappa$ and the commutators of the components of $a$.

**Theorem 1.1**: *There exists $\varepsilon > 0$ (depending at most on the integer n and the vector space $\mathbb{V}$) with the following significance: Let $(A, a)$ denote a pair of connection on $\mathbb{R}^n \times SU(2)$ and section of the associated bundle with fiber $\mathbb{V} \otimes \mathfrak{su}(2)$) obeying (1.1). Either $a \wedge a = 0$ or $\liminf_{r \to \infty} \frac{1}{r^\varepsilon} \kappa(r) > 0$.*

The condition $a \wedge a = 0$ says that $[a_a, a_b] = 0$ for all indices a, b. This implies that $a$ can be written at any point where it is not zero as $a = v\sigma$ with $v \in \mathbb{V}$ and $\sigma \in \mathfrak{su}(2)$.

By way of examples, equations such as those in (1.1) when n = 4 are obeyed by pairs (A, $a$) that obey any of the 1-parameter family of Kapustin-Witten equations [KW]. The version of $\mathbb{V}$ in this case is $\mathbb{R}^4$; and $\mathbb{R}^4 \times \mathbb{V}$ is viewed as $T^*\mathbb{R}^4$ so as to view $a$ as an $\mathfrak{su}(2)$ valued 1-form. The $\tau \in [0, 1]$ version of the Kapustin-Witten equations requires that (A, $a$) obey

- $(1-\tau)(F_A - a \wedge a)^+ - \tau(d_A a)^+ = 0$,
- $\tau(F_A - a \wedge a)^- + (1-\tau)(d_A a)^- = 0$,
- $d_A *a = 0$.

(1.3)

where $F_A$ denotes the curvature 2-form of A, $d_A$ denotes the exterior covariant derivative, and the superscripts $(\cdot)^\pm$ denote the respective self-dual $(\cdot)^+$ and anti-self-dual $(\cdot)^-$ parts of the indicated $\mathfrak{su}(2)$-valued 2-form. What is denoted by $*$ here (and in the upcoming (1.4)) is the Euclidean metric's Hodge star. See, e.g. [T1] for how to get (1.1) from (1.3).

The equation in (1.1) is also obeyed when (A, $a$) are solutions of the Vafa-Witten equations [VW]. These are equivalent on $\mathbb{R}^4$ to the $\tau = 0$ version of (1.3) but they are usually depicted by writing $\mathbb{V}$ as $\mathbb{R}^3 \oplus \mathbb{R}$ with $\mathbb{R}^3$ viewed as the fiber of the self-dual subbundle in $\Lambda^2 T^* \mathbb{R}^4$. Writing $a$ accordingly as $(\mathfrak{a}, \phi)$, these equations say that

- $(F_A)^+{}_a = \frac{1}{2\sqrt{2}} \varepsilon^{abc} [\mathfrak{a}_b, \mathfrak{a}_c] + \frac{1}{\sqrt{2}} [\phi, \mathfrak{a}_a]$,
- $*d_A \mathfrak{a} + d_A \phi = 0$,

(1.4)

with the notation as follows: The indices on $(F_A)^+$ and $\mathfrak{a}$ (which is the self-dual 2-form part of $a$) denote their components with respect to an oriented, orthonormal frame for the self-dual subbundle in $\wedge^2 T^* \mathbb{R}^4$. Meanwhile, $\{\varepsilon^{abc}\}_{a,b,c=1,2,3}$ are the components of the completely anti-symmetric 3-tensor with $\varepsilon^{123} = 1$.

The following theorem says more about the case when n = 4 and (A, $a$) in (1.1) obeys the $\tau = 0$ of (1.3) (which is equivalent to (1.4)) or the $\tau = 1$ version of (1.3).

**Theorem 1.2**: *There exists $\varepsilon > 0$ with the following significance: Suppose that (A, $a$) is a solution on $\mathbb{R}^4$ to either the $\tau = 0$ or $\tau = 1$ version of (1.3) (or to (1.4)). Then either $a \wedge a = 0$ and $\nabla_A a = 0$ or $\lim\text{-}\inf_{r \to \infty} \frac{1}{r^\varepsilon} K(r) > 0$.*

The condition $\nabla_A a = 0$ says that $a$ is A-covariantly constant. In the case when (A, $a$) obeys (1.3), the conditions $a \wedge a = 0$ and $\nabla_A a = 0$ imply that $F_A$ has zero self-dual part when $\tau = 0$ and that it has zero anti-self dual part whent $\tau = 1$. In either case, the curvature need not be identically zero. (If the curvature is non-zero, then it can be written as $F_A = \omega \sigma$ with $\sigma$ being an A-covariantly constant, norm 1 map to $\mathfrak{su}(2)$ and with $\omega$ being an $\mathbb{R}$-valued, closed anti-self dual (or self-dual) 2-form.)

The next theorem talks about the case when $(A, a)$ obeys (1.3) with $\tau \notin \{0, 1\}$.

**Theorem 1.3**: *For any $\tau \in (0, 1)$, there are solutions to the corresponding version of (1.3) with*
a) $\lim_{r \to \infty} K(r) = \sqrt{2}\pi.$
b) $a \wedge a = 0$.
c) $\nabla_A a$ *is not identically zero; its $L^2$ norm on the ball of radius r centered at the origin is greater than a non-zero multiple of r when r is large.*
*With regards to Item c), there are no non-trivial solutions to $\tau \in (0, 1)$ versions (1.3) with $a \wedge a = 0$, with square integrable curvature 2-form, and with $\liminf_{r \to \infty} \frac{1}{r} K(r) = 0$.*

The remainder of this paper is divided into sections that are labeled 2-5. Sections 2-4 prove Theorem 1.1 and Section 5 proves Theorems 1.2 and 1.3.

With regards to notation, the convention in what follows is to use $c_0$ to denote a number that is greater than 1 and independent of any given version of $(A, a)$. It is allowed to depend on the vector space $\mathbb{V}$ and on the dimension n. The precise value of $c_0$ can be assumed to increase between successive appearances.

A second convention concerns the notation for the trace on the vector space of (complex) $2 \times 2$ matrices: Supposing for the moment that $\mu$ is a given matrix, then $\langle \mu \rangle$ denotes $-\frac{1}{2} \text{trace}(\mu)$. This definition is designed so that the bilinear form $(\mathfrak{b}, \mathfrak{c}) \to \langle \mathfrak{b}\mathfrak{c} \rangle$ on the lie algebra $\mathfrak{su}(2)$ is positive definite and thus defines an inner pruduct. Norms are defined by this inner product. With regards to the notation for inner products: If $v$ and $w$ are given elements of $\mathbb{V}$, or tensors on $\mathbb{R}^n$, or $\mathbb{V}$-valued tensors, or $\mathbb{V} \otimes \mathfrak{su}(2)$-valued tensors, then $\langle v, w \rangle$ is used to denote their inner product as defined using the inner product on $\mathbb{V}$, and/or the Euclidean inner product on tensors, and/or the inner product defined just now on $\mathfrak{su}(2)$.

## 2. The function N and the matrix T

The various subsections that follow introduce certain functions on $(0, \infty)$ (viewed as the radius r from the origin in $\mathbb{R}^n$) that are defined using the section $a$. These functions play a central role in the proof of Theorem 1.1.

### a) The frequency function

The derivative of the function $K$ on $(0, \infty)$ where it is non-zero can be written (with the help of (1.1)) as

$$\frac{d}{dr} K = \frac{N}{r} K$$

(2.1)

with $N$ being the function defined where $K > 0$ by the rule

$$N = \frac{1}{r^{n-2}K^2} \int_{B_r} (|\nabla_A a|^2 + \langle [a_a, a], [a_a, a] \rangle) .$$

(2.2)

The function $N$ is called the *frequency* function by virtue of it being the analog of a function with this name and a similar definition that was introduced by Almgren in [A]. The formula in (2.2) indicates that $N \geq 0$, so $K$ is non-decreasing. In particular, if $K$ is ever positive, then it stays positive. A version of Aronszajn's unique continuation theorem is stated and proved in Section 2c which says that $K$ is either identically zero or strictly positive on $(0, \infty)$. Take this on faith for now and assume henceforth that $K > 0$ where $r$ is positive. Supposing that $r > s > 0$, then integrating (2.1) leads to the formula

$$K(r) = \exp\left(\int_s^r \frac{N(s)}{s} ds\right) K(s) .$$

(2.3)

This implies, in particular that if $N(\cdot) \geq \varepsilon$ on $[s, r]$, then $K(r) \geq \left(\frac{r}{s}\right)^\varepsilon K(s)$.

The following formula for the derivative of $N$ is also need:

$$\frac{d}{dr} N = \frac{1}{r^{n-2}K^2} \int_{\partial B_r} (|\nabla_A a|^2 + \langle [a_a, a], [a_a, a] \rangle) - \frac{1}{r}(n - 2 + 2N)N .$$

(2.4)

This formula implies in particular that if $N$ is ever small (say $N \leq \frac{1}{100}$) at some $r > 0$, then it can't be very much larger at points $s < r$ until $\frac{s}{r}$ is relatively small. To be precise:

$$N(s) \leq \left(\frac{r}{s}\right)^n N(r) \quad \text{if} \quad N(r) < 1 \quad \text{and} \quad N(r)^{1/n} r \leq s \leq r.$$

(2.5)

The next lemma makes use of these formulas:

**Lemma 2.1**: *There exists $\kappa > 1$ with the following significance: Fix $\varepsilon \in (0, \kappa^{-1})$. Suppose that $(A, a)$ is a solution to (1.1) and $\rho > \kappa$ is such that the corresponding version of the function $K$ obeys $K(\rho) \leq \rho^\varepsilon K(1)$. There exists $r \in [\rho^{(1-2\sqrt{\varepsilon})}, \rho]$ such that*

- $N(r) \le \sqrt{\varepsilon}$.
- $N < \varepsilon^{1/4}$ *on all of the interval* $[\varepsilon^{1/4n} r, r]$.
- $K \ge (1 - \kappa \varepsilon^{1/4} |\ln \varepsilon|) K(r)$ *on all of the interval* $[\varepsilon^{1/4n} r, r]$.

***Proof of Lemma 2.1***: If r exists with $N(r) < \sqrt{\varepsilon}$ (which is the first bullet's assertion), then (2.5) holds supposing that s is less than r but greater than $\varepsilon^{1/(2n)} r$. This implies in turn that $N(\cdot) < \varepsilon^{1/4}$ on $[\varepsilon^{1/(4n)} r, r]$ which is the second bullet's assertion. The third bullet follows from the second bullet and the identity in (2.3).

With the preceding understood, Lemma 2.1 follows with a proof there is a value of r between $\rho^{(1-2\sqrt{\varepsilon})}$ and $\rho$ with $N(r) \ge \sqrt{\varepsilon}$. To this end, suppose that N is greater than $\sqrt{\varepsilon}$ on the whole of the interval between $\rho^{1-2\sqrt{\varepsilon}}$ and $\rho$. If this is so, then the right hand side of (2.3) is greater than $\rho^{2\varepsilon}$ which violates the assumption that $K(\rho) < \rho^{\varepsilon}$. □

Supposing that $v$ is a constant, unit length vector in $\mathbb{V}$, let $a(v)$ denote $a_c v_c$. This is an $\mathfrak{su}(2)$-valued function on $\mathbb{R}^n$ that obeys

$$\nabla_A^{\dagger} \nabla_A a(v) + [a_c, [a(v), a_c]] = 0 .$$

(2.6)

The function $a(v)$ has its own version of K to be denoted by $K_v$ whose square is given by the formula in (1.2) with $a$ replaced by $a(v)$. Section 2c proves a unique continuation theorem for $a(v)$ also which implies that either $a(v)$ is identically zero or $K_v > 0$ on $(0, \infty)$. Supposing that $a(v)$ is not identically zero, then the derivative of the function $K_v$ can be written like that of K where $K_v > 0$:

$$\frac{d}{dr} K_v = \frac{N_v}{r} K_v$$

(2.7)

with $N_v$ defined by the formula

$$N_v = \frac{1}{r^{n-2} K_v^2} \int_{B_r} (|\nabla_A a(v)|^2 + \|[a, a(v)]\|^2) .$$

(2.8)

Likewise, (2.3) holds with $K_v$ replacing K and $N_v$ replacing N. And, as was the case with N, if the function $N_v$ is small at some value of r (say $N_v(r) < 1$), then it can not be too big at $s < r$. In particular, (2.5) holds with $N_v$ replacing N on both sides. (This $N_v$ version of (2.5) is proved by differentiating the formula in (2.8) for $N_v$ to see that its derivative is no smaller than $-\frac{n}{r} N_v$ where $N_v < 1$.)

**b) The matrix $\mathbb{T}$**

What is denoted below by $\mathbb{T}$ is a $\mathbb{V} \otimes \mathbb{V}$ valued function on $[0, \infty)$ whose components with respect to an orthonormal basis of $\mathbb{V}$ are defined by the rule:

$$\mathbb{T}_{ab}(r) = \frac{1}{r^3} \int_{\partial B_r} \langle a_a a_b \rangle .$$

(2.9)

This matrix is symmetric and non-negative definite. Its trace is $\kappa^2$ and, if $v$ is any given unit length vector in $\mathbb{V}$, then $\mathbb{T}_{ab} v_a v_b = \kappa_v^2$. Define the norm of $\mathbb{T}$ (to be denoted by $|\mathbb{T}|$) by setting its square to be $\mathbb{T}_{ab} \mathbb{T}_{ab}$. This norm obeys $\frac{1}{3}\kappa^2 \le |\mathbb{T}| \le \kappa^2$. The matrix $\mathbb{T}$ is differentiable on $(0, \infty)$ and the norm of its derivative obeys

$$|\tfrac{d}{dr} \mathbb{T}| \le c_0 \tfrac{N}{r} |\mathbb{T}|.$$

(2.10)

No generality is lost by assuming that $\mathbb{T}$ does not have a zero eigenvalue at positive r. To explain: If $\mathbb{T}$ has a zero eigenvalue at some positive radius r, then the corresponding eigenvector (call it $v$) is such that $\kappa_v^2(r) = 0$. As a consequence of the unique continuation asserted in Section 2c, this same $v$ has $\kappa_v(\cdot)$ being identically zero and thus $a(v)$ identically zero. In this case, $\mathbb{V}$ can be replaced by the orthogonal complement to $v$ in $\mathbb{V}$ and (1.1) holds for this smaller vector space. In this way, $\mathbb{V}$ can be cut down in size so that the resulting version of $\mathbb{T}$ has no zero eigenvalues.

Given $r \in (0, \infty)$, let $\lambda(r)$ denote the smallest eigenvalue of $\mathbb{T}(r)$. The following lemma says more about $\lambda$:

**Lemma 2.2**: *Let $(A, a)$ denote a solution to (1.1) with the property that the matrix $\mathbb{T}$ has no zero eigenvalue at positive* r. *Let $\lambda(\cdot)$ denote the function on $(0, r)$ whose value at any given* r *is the smallest eigenvalue of $\mathbb{T}(r)$. This $\lambda$ is a Lipschitz function on $(0, \infty)$. Moreover, $\lambda$ is nearly differentiable in the following sense: Given $r \in (0, \infty)$, let $v \in \mathbb{V}$ denote a unit length eigenvector of $\mathbb{T}(r)$ with the eigenvalue $\lambda(r)$. Then*

- $\lambda(r+\Delta) - \lambda(r) \le \langle v, \mathbb{T}(r+\Delta) v \rangle - \langle v, \mathbb{T}(r) v \rangle = \langle v, (\tfrac{d}{dr} \mathbb{T})_r v \rangle \Delta + \mathcal{O}(\Delta^2)$.
- $\lambda(r) - \lambda(r-\Delta) \ge \langle v, \mathbb{T}(r) v \rangle - \langle v, \mathbb{T}(r-\Delta) v \rangle = \langle v, (\tfrac{d}{dr} \mathbb{T})_r v \rangle \Delta + \mathcal{O}(\Delta^2)$.

*Proof of Lemma 2.2*: Given r and $\Delta$, let $v$ denote an eigenvector of $\mathbb{T}(r)$ with eigenvalue $\lambda(r)$ and let $v_\Delta$ denote an eigenvector of $\mathbb{T}(r+\Delta)$ with eigenvalue $\lambda(r+\Delta)$. Then

$$\lambda(r+\Delta) - \lambda(r) = \langle v_\Delta, \mathbb{T}(r+\Delta) v_\Delta \rangle - \langle v, \mathbb{T}(r) v \rangle$$

(2.11)

by definition. However, because $\lambda(r+\Delta)$ is the smallest eigenvalue of $\mathbb{T}$, the right hand side of (2.11) is no smaller than $\langle v, \mathbb{T}(r+\Delta)v \rangle - \langle v, \mathbb{T}v \rangle$. This is the left most inequality in the top bullet of the lemma. The right most equality follows via Taylor's theorem with remainder. An $\mathcal{O}(\Delta)$ lower bound for $\lambda(r+\Delta) - \lambda(r)$ (and this a proof that $\lambda$ is Lipschitz) follows from the identity in (2.11) by virtue of the fact that $\langle v, \mathbb{T}(r) v \rangle$ is no larger than $\langle v_\Delta, \mathbb{T}(r) v_\Delta \rangle$. This implies that $\lambda(r+\Delta) - \lambda(r)$ is no less than $\langle v_\Delta, (\mathbb{T}(r+\Delta) - \mathbb{T}(r)) v_\Delta \rangle$ whose norm is at most $c_0 \Delta |\frac{d}{dr} \mathbb{T}| + \mathcal{O}(\Delta^2)$. The argument for the second bullet of the lemma is much like the argument for the first. □

This lemma has the following implication: Suppose that $r > s$ are two positive numbers. Let N denote a positive integer. For each $k \in \{0, 1, \ldots, N\}$, define the number $r_k$ by the rule $r_k = s + \frac{k}{N}(r-s)$. For each such k, let $v_k$ denote an eigenvector of $\mathbb{T}(r_k)$ with eigenvalue $\lambda(r_k)$. The two bullets in Lemma 2.2 imply that $\lambda(r)$ can be written as:

$$\lambda(r) = \lambda(s) + \sum_{k=1}^{N} \langle v_k, (\tfrac{d}{dr}\mathbb{T})_{r_k} v_k \rangle \tfrac{r-s}{N} + \mathcal{O}(\tfrac{1}{N}).$$

(2.12)

Meanwhile, $\langle v_k, (\frac{d}{dr}\mathbb{T})_{r_k} v_k \rangle$ can be written (using (2.7) and (2.8)) as

$$\langle v_k, (\tfrac{d}{dr}\mathbb{T})_{r_k} v_k \rangle = \tfrac{1}{r_k} \lambda(r_k) N_{v_k}(r_k)$$

(2.13)

with $N_{v_k}$ being the $v = v_k$ version of the function $N_v$. The equation in (2.12) can be written using this notation as

$$\lambda(r) = \lambda(s) + \sum_{k=1}^{N} \tfrac{1}{r_k} \lambda(r_k) N_{v_k}(r_k) \tfrac{r-s}{N} + \mathcal{O}(\tfrac{1}{N}).$$

(2.14)

Since $N_v$ for any $v \in \mathbb{V}$ is nonnegative, the equation in (2.14) implies that $\lambda$ is a non-decreasing function of r.

The next lemma adds to a slightly weakened version of Lemma 2.1.

**Lemma 2.3**: *There exists $\kappa > 100$ with the following significance: Supposing that $(A, a)$ is a solution to (1.1), use a to define the function $\kappa$, the matrix $\mathbb{T}$ and the function $\lambda$. Fix $\varepsilon \in (0, \kappa^{-1})$ and suppose that $\rho$ is greater than $\kappa(1 + \lambda(1)^{-2/\varepsilon})$, and that the value of $\kappa$ at $\rho$ obeys $\kappa(\rho) \leq \rho^\varepsilon \kappa(1)$. Then, there exists $r \in [\rho^{(1-30\sqrt{\varepsilon})}, \rho]$ such that*

- *The functions $N$ and $K$ obey*
  a) $N < \varepsilon^{1/4}$ on all of the interval $[\varepsilon^{1/(8n)} r, r]$.
  b) $K \geq (1 - \kappa \varepsilon^{1/4} |\ln \varepsilon|) K(r)$ on all of the interval $[\varepsilon^{1/(8n)} r, r]$.
- *Let $v$ denote an eigenvector of $\mathbb{T}(r)$ with the smallest of $\mathbb{T}$'s eigenvalues. The corresponding functions $N_v$ and $K_v$ obey*
  a) $N_v \leq \varepsilon^{1/4}$ on all of $[\varepsilon^{1/(8n)} r, r]$.
  b) $K_v \geq (1 - \kappa \varepsilon^{1/4} |\ln \varepsilon|) K_v(r)$ on all of the interval $[\varepsilon^{1/(8n)} r, r]$.

*Proof of Lemma 2.3*: The proof has four steps.

Step 1: This step makes an observation that was implicit in Lemma 2.1's proof. To state this observation, introduce by way of notation $\Omega$ to denote the set $x \in [0, \ln \rho]$ where $N(e^x) > \sqrt{\varepsilon}$. The measure of $\Omega$ can not be greater than $2\sqrt{\varepsilon} \ln \rho$ because, otherwise, (2.3) would have $\ln(K(\rho)/K(1))$ being greater than $2\varepsilon \ln \rho$ which violates the assumptions.

Step 2: Consider (2.14) in the case when $r = 2s$. Since $\lambda$ is a non-decreasing function of $r$, the $r = 2s$ version of (2.14) leads to the inequality:

$$\frac{\lambda(2s)}{\lambda(s)} \geq \left(1 + \frac{1}{2N} \sum_{k=1}^{N} N_{v_k}(r_k) + \mathcal{O}\left(\frac{1}{N}\right)\right).$$

(2.15)

Now let L denote the largest integer such that $2^L < \rho$. Taking $s$ to be successively $2^m$ for $m = 0, 1, \ldots, L-1$ and invoking (2.15) in each case leads to the inequality

$$\frac{\lambda(\rho)}{\lambda(1)} \geq \prod_{m=0}^{L-1} (1 + x_m)$$

(2.16)

where $x_m$ shorthand for $x_m = \frac{1}{2N} \sum_{k=1}^{N} N_{v_k}(2^m (1 + \frac{k}{N})) + \mathcal{O}(\frac{1}{N})$.

Step 3: Let $f$ denote the fraction of integers in the set $\{0, 1, \ldots, L-1\}$ where the corresponding version of $x_m$ is greater than $\sqrt{\varepsilon}$. Using just this fraction, (2.16) leads to

$$\frac{\lambda(\rho)}{\lambda(1)} \geq (1 + \sqrt{\varepsilon})^{fL}.$$

(2.17)

Meanwhile, $\lambda(\rho) < \rho^\varepsilon$ so (2.17) implies (by taking logarithms) that

$$\varepsilon \ln \rho - \ln(\lambda(1)) \geq \tfrac{1}{2} \sqrt{\varepsilon} f L$$

(2.18)

if $\varepsilon < \frac{1}{100}$. Since $L > \frac{1}{2\ln 2} \ln \rho$ if $\rho > 8$ so (2.18) gives a bound on $f$:

$$f \leq 4\ln(2)\sqrt{\varepsilon}(1 - \frac{\ln(\lambda(1))}{\ln(\rho^\varepsilon)}) \, .$$

(2.19)

Therefore, if $\rho \geq \lambda(1)^{-2/\varepsilon}$, then $f \leq 8\ln(2)\sqrt{\varepsilon}$.

<u>Step 4</u>: A given integer $m \in \{0, 1, \ldots, L-1\}$ labels a corresponding interval in the set $[0, \ln(\rho)]$, this being $[m\ln(2), (m+1)\ln(2)]$. This understood, it follows from what is said in Step 3 that the subset of $[0, \ln(\rho)]$ accounted for by integers $m \in \{1, 2, \ldots, L-1\}$ with $x_{m-1} > \sqrt{\varepsilon}$ has measure at most $8\sqrt{\varepsilon}\ln(\rho)$. Let $\Omega_\lambda$ denote this subset of $[0, \ln(\rho)]$.

The union of $\Omega$ (from Step 1) and $\Omega_\lambda$ has measure at most $10\sqrt{\varepsilon}\ln(\rho)$. As a consequence, there is a point (call it $r'$) not in $\Omega_\lambda$ and not in $\Omega$ and from the interval $[\rho^{1-20\sqrt{\varepsilon}\ln(\rho)}, 2^L]$. Thus, $N(r') < \sqrt{\varepsilon}$. Also, this point $r'$ is in some interval of the form $[2^m, 2^{m+1}]$ with $m$ being an integer from the set $\{1, \ldots, L\}$. Since the corresponding $x_{m-1}$ is less then $\sqrt{\varepsilon}$, there exists some $r$ in the interval $[2^{m-2}, 2^{m-1}]$ and eigenvector $v$ of $\mathbb{T}(r)$ with the smallest eigenvalue (thus $\lambda(r)$) such that $N_v(r) < 2\sqrt{\varepsilon}$. (This $r$ can be taken to be one of the $r_k$'s that appear in the large $N$ version of the formula for $x_m$.) It follows from Lemma 2.1 that $N(r) \leq \varepsilon^{1/4}$ and that $N < \varepsilon^{1/4}$ on the whole of the interval $[\varepsilon^{1/8n} r, r]$. Meanwhile, it follows from (2.3) that $K \geq (1 - c_0 \varepsilon^{1/4} |\ln \varepsilon|) K(r)$ on this same interval. The $N_v$ version of (2.5) and the identity in (2.7) can be used to obtain the analogous bounds for $N_v$ and $K_v$.

**c) Unique continuation**

The purpose of this subsection is to tie up a loose end by proving that any $v \in \mathbb{V}$ version of the function $K_v$ is either strictly positive on $(0, \infty)$ or identically zero. This will imply the same dichotomy for the function $K$.

To this end, note that the set where $K_v > 0$ has the form $(r_0, \infty)$ since $K_v$ is continuous. If $N_v$ is a priori bounded on $(r_0, 2r_0]$, then (2.7) prevents $K_v$ from vanishing at $r_0$ if $r_0$ is positive. Now, $N_v$ is, in any event, bounded by some number (call it $\mathcal{E}$) at $2r_0$ so it is sufficient to prove that $N_v$ can not be much greater than $\mathcal{E}$ on $(r_0, 2r_0]$ if $r_0 > 0$. The proof that this is the case starts with the $N_v$ version of (2.4) which is written explicitly below:

$$\frac{d}{dr} N_v = \frac{1}{r^{n-2} K_v^2} \int_{\partial B_r} (|\nabla_A a(v)|^2 + |[a, a(v)]|^2) - \frac{1}{r}(n - 2 + 2N_v) N_v \, .$$

(2.20)

This inequality by itself is not sufficient to prove that $N_v$ does not diverge as $r \to r_0$. A rewriting of (2.20) is needed for this purpose. The rewriting of (2.20) and then its application to bounding $N_v$ has four parts. This rewriting and the subsequent bounding of $N_v$ copies for the most part what is done by Almgren [A].

*Part 1*: Fix Euclidean coordinates $\{x_\alpha\}_{\alpha=1,\ldots,n}$ for $\mathbb{R}^n$. Having done so, introduce a symmetric, $n \times n$ tensor (to be denoted by S) by the rule that sets its components to be

$$S_{\alpha\beta} = \langle \nabla_{A\alpha} a(v)_c, \nabla_{A\beta} a(v)_c \rangle - \tfrac{1}{2} \delta_{\alpha\beta} (|\nabla_A a(v)|^2 + |[a, a(v)]|^2)$$

(2.21)

with the notation as follows: The collection $\{\nabla_{A\alpha}\}_{\alpha=1,\ldots,n}$ are the directional covariant derivatives in the coordinate axis directions; and $\delta_{\alpha\beta}$ is 1 if $\alpha \neq \beta$ and 0 otherwise. This tensor would be divergence free were it not for the curvature of A and the other components of $a$: A computation using (2.6) writes the divergence of S as

$$\nabla_\beta S_{\alpha\beta} = \langle F_{A\beta\alpha}, [a(v)_c, \nabla_{A\beta} a(v)_c] \rangle - \langle [\nabla_{A\alpha} a_c, a(v)][a_c, a(v)] \rangle$$

(2.22)

with $\{F_{A\beta\alpha}\}_{\beta,\alpha \in \{1,2,3,4\}}$ denoting the components of $F_A$.

*Part 2*: View (2.22) as an equality of two vectors fields on $\mathbb{R}^n$. Take the inner product of both of these vector fields with the radial pointing vector $x^\alpha \frac{\partial}{\partial x^\alpha}$ and then integrate the resulting equality of functions over the ball $B_r$. An integration by parts writes the integral identity as

$$r \int_{\partial B_r} (|\nabla_{Ar} a(v)|^2 - \tfrac{1}{2} |\nabla_A a(v)|^2 - \tfrac{1}{2} |[a, a(v)]|^2) + \int_{B_r} (\tfrac{n-2}{2} |\nabla_A a|^2 + \tfrac{n}{2} |[a, a(v)]|^2) = \int_{B_r} x^\alpha \langle F_{A\beta\alpha}, [a(v)_c, \nabla_{A\beta} a(v)_c] \rangle + \int_{B_r} x^\alpha \langle [\nabla_{A\beta} a_c, a(v)], [a_c, a(v)] \rangle .$$

(2.23)

Here, $\nabla_{Ar}$ denotes the directional covariant derivative along the unit vector tangent to the lines through the origin. The terms in this identity can be rearranged to say that:

$$\int_{\partial B_r} (|\nabla_A a(v)|^2 + |[a, a(v)]|^2) = 2 \int_{\partial B_r} |\nabla_{Ar} a(v)|^2 + \frac{(n-2)}{r} \int_{B_r} (|\nabla_A a(v)|^2 + \tfrac{n}{n-2} |[a, a(v)]|^2) - \tfrac{2}{r} \mathcal{E}$$

(2.24)

with $\mathcal{E}$ denoting here the two integrals that appear on the right hand side of (2.23). Use this identity for the boundary integral on the right hand side of (2.20). What with the definition of $N_v$ in (2.8), this rewriting leads in turn to an inequality asserting that

$$\frac{d}{dr} N_v \geq \frac{1}{r^{n-2}K_v^2} \int_{\partial B_r} (|\nabla_{Ar}a(v) - \frac{N_v}{r}a(v)|^2 - \frac{2}{r^{n-1}K_v^2}\mathfrak{E} .$$

(2.25)

If the $\mathfrak{E}$ term were absent, then this inequality says that $N_v$ is an increasing function of r and so if it is bounded at $2r_0$, then it is bounded at $r_0$ also. This conclusion also follows with a suitable bound on $|\mathfrak{E}|$ which is what is done in Part 3.

*Part 3*: The $\mathfrak{E}$ term in (2.25) is a sum of two integrals, the first being

$$\int_{B_r} x^\alpha \langle F_{A\beta\alpha}, [a(v)_c, \nabla_{A\beta}a(v)_c] \rangle$$

(2.26)

Supposing that $r < 2r_0$, then there is a bound on $|F_A|$ on $B_r$ by some r-independent number to be denoted by $\mathcal{F}$. The size of $\mathcal{F}$ is of no consequence. Use of this bound leads to the following bound on the absolute value of integral in (2.26):

$$|\int_{B_r} x^\alpha \langle F_{A\beta\alpha}, [a(v)_c, \nabla_{A\beta}a(v)_c] \rangle| \leq c_0 \, r \mathcal{F} (\int_{B_r} |a(v)|^2)^{1/2} (\int_{B_r} |\nabla_A a(v)|^2)^{1/2} .$$

(2.27)

Since $K_v$ is a non-decreasing function, the $L^2$ norm of $a(v)$ that appears on the right hand side of (2.27) is no greater than $r^{n/2}K_v(r)$. Meanwhile, the $L^2$ norm of $\nabla_A a(v)$ that appears is at most $r^{(n-2)/2}K_v\sqrt{N_v}$. Hence, the right hand side of (2.27) is at most $c_0 \mathcal{F} r^n K_v^2 \sqrt{N_v}$.

The second integral that appears on the right hand side of (2.23) is

$$\int_{B_r} x^\alpha \langle [\nabla_{A\beta}a_c, a(v)], [a_c, a(v)] \rangle .$$

(2.28)

Both $|a|$ and $|\nabla_A a|$ have some r-independent upper bound on $B_r$ if $r \leq 2r_0$. This can be assumed to be $\mathcal{F}^{1/2}$ (at the expense of making $\mathcal{F}$ bigger). This being the case, the integral in (2.28) is at most $c_0 \mathcal{F} r^{n+1} K_v^2(r)$.

*Part 4*: Using the bounds from Part 3 in (2.25) leads to the inequality

$$\frac{d}{dr} N_v \geq -c_0 \mathcal{F} (r^2 + r\sqrt{N_v}) .$$

(2.29)

This in turn says that the derivative of $N_v$ is no smaller than $-c_0\mathcal{F}(r^2 + N_v)$ and this lower bound implies that $N_v(r)$ can not diverge as $r \to r_0$ if it is finite at $r = 2r_0$. And, given that $N_v$ is bounded on $(r_0, 2r_0)$ and that $K_v > 0$ at $r = 2r_0$, then (2.7) won't let $K_v$ vanish at $r = r_0$ when $r_0$ is positive.

### 3. Pointwise bounds on $a$

Let $v$ denote a chosen, unit norm vector in $\mathbb{V}$. The proposition that follows asserts a priori pointwise bounds for $a(v)$ on balls about the origin in $\mathbb{R}^n$ given suitable bounds on $K_v$ and $N_v$. By way of notation, the proposition uses $\omega$ to denote the volume of the radius 1 sphere about the origin in $\mathbb{R}^n$.

**Proposition 3.1**: *There exists $\kappa > 100$ with the following significance: Suppose that $(A, a)$ is a solution to (1.1). Let $v \in \mathbb{V}$ denote a unit length vector. Fix $r > 0$.*

- $|a(v)| \leq \kappa K_v(r)$ *on the ball of radius $\frac{7}{8} r$ centered at the origin.*
- $|a(v)| \leq \frac{1}{\sqrt{\omega}} (1 + \kappa \sqrt{N_v(r)}) K_v(r)$ *on the ball of radius $\frac{7}{8} r$ centered at the origin.*
- *Supposing that $s \in [\frac{1}{2} r, \frac{7}{8} r]$, the set of points in the radius $s$ sphere about the origin where $|a(v)| \leq \frac{1}{2\sqrt{\omega}} K_v(r)$ has $(n-1)$ dimensional volume at most $\kappa \sqrt{N_v(r)} \, r^{n-1}$.*

By way of a parenthetical remark: The statements of the second and third bullets of the lemma are useful only in the case when $N_v(r) \leq c_0^{-1}$. By way of a second remark: The bounds asserted by the three bullets of the proposition hold on any concentric ball inside $B_r$ with strictly smaller radius with it understood that the number $\kappa$ depends on the number $\mu$ from $(0, 1)$ that is defined by writing the radius of this nested ball as $(1 - \mu) r$.

As is explained in Section 3d, the full norm of $a$ obeys bounds much like the bounds from Proposition 3.1 for its components. The following proposition makes a precise statement to this effect. (This proposition is not used in subsequent arguments.)

**Proposition 3.2**: *There exists $\kappa > 100$ with the following significance: Suppose that $(A, a)$ is a solution to (1.1). Fix $r > 0$.*

- $|a| \leq \kappa K(r)$ *on the ball of radius $\frac{7}{8} r$ centered at the origin.*
- $|a| \leq \frac{1}{\sqrt{\omega}} (1 + \kappa \sqrt{N(r)}) K(r)$ *on the ball of radius $\frac{7}{8} r$ centered at the origin.*

*Supposing that $s \in [\frac{1}{2} r, \frac{7}{8} r]$, the set of points in the radius $s$ sphere about the origin where $|a| \leq \frac{1}{2\sqrt{\omega}} K(r)$ has $(n-1)$ dimensional volume at most $\kappa \sqrt{N(r)} \, r^{n-1}$.*

The proof of Proposition 3.1 is modeled on a proof of similar assertions in [T1] and [T2]. The proof occupies Sections 3a-c. Section 3d has the proof of Proposition 3.2 (which also has close kin in these same references).

### a) The first a priori pointwise bound on $|a(v)|$

This subsection proves the top bullet of Proposition 3.1. The proof has four steps.

<u>Step 1</u>: Fix a smooth, non-increasing function on $\mathbb{R}$ to be denoted by $\beta$ that is equal to 1 on $(-\infty, \frac{15}{16}]$ and equal to 0 on $[\frac{31}{32}, \infty)$. Supposing that $r > 0$, let $\beta_r$ denote the function on $\mathbb{R}^n$ given by the rule $x \to \beta(\frac{1}{r}|x|)$. This function is equal to 1 where $|x|$ is less than $\frac{15}{16}r$ and it is equal to 0 where $|x|$ is greater than $\frac{31}{32}r$. In particular, $\beta_r$ has compact support in $B_r$. Notice also that $|d\beta_r| \leq c_0 \frac{1}{r}$.

<u>Step 2</u>: Given $p \in B_r$, let $G_p$ denote the Dirichelet Green's function for the Laplacian (this is $d^\dagger d$) on $B_r$ with pole at the point $p$. This is a smooth, non-negative function on $B_r - p$. If the dimension n is greater than 2, then $G_p$ obeys

$$G_p(x) \leq \frac{1}{\omega} \frac{1}{|x-p|^{n-2}} \quad \text{and} \quad |dG_p| \leq c_0 \frac{1}{|x-p|^{n-1}} \ .$$

(3.1)

If the dimension $n = 2$, then $G_p(x) \leq c_0 \ln(\frac{r}{|x-p|})$ and $|dG_p| \leq c_0 \frac{1}{|x-p|}$.

<u>Step 3</u>: Take the inner product of both sides of (2.6) with $a(v)$ to obtain

$$\tfrac{1}{2} d^\dagger d |a(v)|^2 + |\nabla_A a(v)|^2 + |[a, a(v)]|^2 = 0$$

(3.2)

Multiply both sides of this identity by $\beta_r G_p$ and then integrate the result over $B_r$. An integration by parts leads from the integral identity to the following:

$$\tfrac{1}{2}|a(v)|^2(p) + \int_{B_r} G_p(|\nabla_A a(v)|^2 + |[a, a(v)]|^2) = -\int_{B_r} (\langle dG_p, d\beta_r \rangle - \tfrac{1}{2} G_p d^\dagger d\beta_r)|a(v)|^2 \ .$$

(3.3)

This identity bounds $|a(v)|$ at $p$ given a suitable bound for the norm of its right hand side.

<u>Step 4</u>: Suppose now that $|p| < \frac{15}{16}r$ so that the distance from p to the support of $d\beta_r$ and $d^\dagger d\beta_r$ is greater than $\frac{1}{16}r$. The right hand side of (3.3) is therefore no greater than $c_0 \frac{1}{r^n}$ times the integral of $|a(v)|^2$ over $B_r$. And, the latter integral is no greater than $c_0 r^n \kappa_v^2(r)$ since $\kappa_v$ is an increasing function. Therefore, the right hand side of (3.3) (and thus $|a(v)|^2$ at p) is not larger than $c_0 \kappa_v^2(r)$.

**b) The second a priori pointwise bound on $|a(v)|$**

This step proves the assertion made by the second bullet of Proposition 3.1. The proof of this bullet has three steps.

<u>Step 1</u>: Define a non-negative function to be denoted by $M_v$ on $\mathbb{R}^4 \times (0, \infty)$ by writing its square at any given pair (p, s) to be

$$M_v^2(p,s) = \frac{1}{s^n} \int_{|(\cdot)-p| \leq s} |a(v)|^2$$

(3.4)

Up to a constant factor, the square of $M_v(p,s)$ is the average of $|a(v)|^2$ on the radius s ball centered at p.

Now, the point p has its own version of the function $K_v$ (to be denoted by $K_v(p;\cdot)$) whose square is defined by the rule

$$K_v^2(p,s) = \frac{1}{s^{n-1}} \int_{|(\cdot)-p|=s} |a(v)|^2 \ .$$

(3.5)

As was the case with the p = 0 version of $K_v$ that is described in Section 2, the point p's version $K_v(p,s)$ is an increasing function on $[0, \infty)$. Since its $s \to 0$ limit is $\sqrt{\omega}\,|a(v)|(p)$, this number is a lower bound for $K_v(p,s)$ for any s. It follows as a consequence that

$$M_v(p,s) \geq \sqrt{\tfrac{\omega}{n}}\,|a(v)|(p) \ .$$

(3.6)

Meanwhile, $\frac{1}{\sqrt{n}} K_v(p,s)$ is an upper bound for $M_v(p;s)$, again because $K_v(p,s)$ is increasing. For example, if r > 0 and if s ≤ r, then $M_v(0,s) \leq \frac{1}{\sqrt{n}} K_v(r)$ with $K_v(r)$ being the p = 0 version of $K_v(p,r)$ (which is the version that appears in Proposition 3.1).

<u>Step 2</u>: Supposing that s is fixed, the resulting function $M_v(\cdot,s)$ on $\mathbb{R}^4$ is differentiable with its derivative obeying

$$M_v|dM_v| \leq \frac{1}{s^n} \int_{|(\cdot)-p| \leq s} |a(v)|\,|\nabla_A a(v)| \ .$$

(3.7)

(Derive this bound by differentiating both sides of the formula in (3.4).) The bound in (3.7) implies in turn that

$$|dM_v| \leq \frac{1}{s^{n/2}} \left( \int_{|(\cdot)-p| \leq s} |\nabla_A a(v)|^2 \right)^{1/2}.$$

(3.8)

If $r > 0$ has been given, and if the point p is in the radius $\frac{7}{8} r$ ball centered at the origin, and if $s \leq \frac{1}{8} r$, then the integral on the right hand side of (3.8) can be compared to the integral that defines the function $N_v$; and doing so leads to the bound

$$|dM_v| \leq \frac{1}{r} \left( \frac{r}{s} \right)^{n/2} K_v(r) \sqrt{N_v(r)}.$$

(3.9)

For example, if $s = \frac{1}{8} r$, then this implies that $|dM_v| \leq c_0 \frac{1}{r} K_v(r) \sqrt{N_v(r)}$.

<u>Step 3</u>: Take $s = \frac{1}{8} r$ and take p to be a point in the radius $\frac{7}{8} r$ ball centered at the origin in $\mathbb{R}^4$. The bound in (3.9) implies that

$$M_v(p, \tfrac{1}{8} r) \leq M_v(0, \tfrac{1}{8} r) + c_0 K_v(r) \sqrt{N_v(r)}.$$

(3.10)

This inequality in (3.10) implies in turn that

$$\sqrt{\tfrac{\omega}{n}} \, |a(v)|(p) \leq \tfrac{1}{\sqrt{n}} (1 + c_0 \sqrt{N_v(r)}) K_v(r)$$

(3.11)

because $M_v(0, \tfrac{1}{8} r) \leq \tfrac{1}{\sqrt{n}} K_v(r)$ (according to Step 2) and because $M_v(p, \tfrac{1}{8} r)$ obeys (3.6). The inequality in (3.11) implies what is asserted by the second bullet of Proposition 3.1.

### c) The set where $|a_v|$ is small

This last part of the subsection proves the third bullet of Proposition 3.1. The proof of this bullet has two steps. The proof assumes at various instances that $N_v(r) \leq c^{-1}$ with $c > 1$ being independent of $(A, a)$. This is sufficient because if $c$ is given a priori (and is independent of $(A, a)$) and if $N_v(r) \geq c^{-1}$, then the third bullet follows directly if the number $\kappa$ is taken to be larger than $\omega c$.

<u>Step 1</u>: It is a consequence of (2.8) that the function $s \to s^{n-2} K_v^2 N_v$ is an increasing function on $(0, \infty)$. Supposing that $r > 0$ has been given; then the preceding observation with (2.7) lead to a bound the derivative of $K_v^2$ at any $s \in (0, r]$ by

$$\left( \tfrac{d}{dr} K_v^2 \right) \big|_s \leq \tfrac{2}{r} \left( \tfrac{r}{s} \right)^{n-1} N_v(r) K_v^2(r).$$

(3.12)

As a consequence, if $s \in [\frac{1}{2} r, \frac{7}{8} r]$, then integrating (3.12) from s to r leads to:

$$K_v^2(s) \geq (1 - c_0 N_v(r)) K_v^2(r)$$

(3.13)

(It is probably needless to say that this inequality is only useful if $N_v(r) \leq c_0^{-1}$.)

Step 2: Introduce $\Omega$ to denote the subset of the radius s sphere about the origin where $|a(v)| \leq \frac{1}{2\sqrt{\omega}} K_v(r)$. Write the n-1 dimensional volume of $\Omega$ as $\omega s^{n-1} f_\Omega$.

The integral that defines $K_v^2(s)$ is the integral of $|a(v)|^2$ over the radius s sphere about the origin. Therefore, it can be bounded from above by replacing the integrand on $\Omega$ by $\frac{1}{4\omega} K_v^2(r)$ and by replacing the integrand on the complement of $\Omega$ by the bound on $|a(v)|^2$ from the second bullet of Proposition 3.1. Doing this leads to the upper bound

$$K_v^2(s) \leq \tfrac{1}{4} K_v^2(r) f_\Omega + (1 + \kappa_* \sqrt{N_v(r)})^2 (1 - f_\Omega) K_v^2(r)$$

(3.14)

with $\kappa_*$ being the version of $\kappa$ that appears in the second bullet of Proposition 3.1. This inequality (after a rearrangment) implies in turn that

$$f_\Omega \leq \tfrac{4}{3} (1 - \frac{1}{(1+\kappa_* \sqrt{N_v(r)})^2} \frac{K_v^2(s)}{K_v^2(r)}) \ .$$

(3.15)

Supposing that $N_v(r) \leq c_0^{-1}$, then (3.13) and (3.15) together lead to the bound that is asserted by the third bullet of Proposition 3.1.

**d) Proof of Proposition 3.2**

The first bullet of Proposition 3.2 follows from the first bullet of Proposition 3.1 since K is no smaller than any given $v \in \mathbb{V}$ version of $K_v$ (assuming $|v| = 1$). The second bullet of Proposition 3.2 follows from the second bullet of Proposition 3.1 by virtue of K being as large as any given version of $K_v$ and by virtue of $NK^2$ being as large as any given version of $N_v K_v^2$ (assuming $|v| = 1$). This last bound follows by directly comparing (2.2) with (2.8). The third bullet follows by virtue of the fact that the set where $|a| \leq \frac{1}{2\sqrt{\omega}} K(r)$ is part of the set where $|a(v)| \leq \frac{1}{2\sqrt{\omega}} K_v(r)$ with v being a norm 1 eigenvector of the matrix $\mathbb{T}(r)$ for its largest eigenvalue. Note in this regard that $K_v(r) \geq c_0^{-1} K(r)$ for this version of v; and this upper bound implies that $N_v(r) \leq c_0 N(r)$.

## 4. Proof of Theorem 1.1

This section contains the proof of Theorem 1.1. To set up the proof, let (A, $a$) denote a solution to (1.1). For any given $r > 0$, use $a$ to define the matrix $\mathbb{T}(r)$ using the rule depicted in (2.9); and assume that this matrix is non-degenerate for all r. (Remember that this can be arranged for any given pair (A, $a$) by restricting to a suitable linear subspace in $\mathbb{V}$.)  Assuming this non-degeneracy condition, and given $r > 0$, let $\lambda(r)$ denote the smallest eigenvalue of $\mathbb{T}(r)$ (which is necessarily a positive number).

With $\kappa_1$ denoting Lemma 2.3's version of the number $\kappa$, suppose that $\varepsilon > 0$, that $\rho > \kappa_1(1 + \lambda(1)^{-2/\varepsilon})$ and that $K(\rho) \leq \rho^\varepsilon K(1)$. The six parts of the proof that follow explain why this assumption leads to nonsense if $\varepsilon < c_0^{-1}$ and if $\rho$ is greater than the $c_0$ times the larger of the numbers $K(1)^{1/(1-30\sqrt{\varepsilon})}$ and $c_0 \kappa_1(1 + \lambda(1)^{-2/\varepsilon})$.

*Part 1*: Let $\kappa_1$ denote the version of the number $\kappa$ from Lemma 2.3. Supposing that $\varepsilon < \kappa_1^{-1}$, then Lemma 2.3 finds a $r \in [\rho^{(1-30\sqrt{\varepsilon})}, \rho]$ so that its two bulleted assertions hold. This value of r is denoted below by $r_1$. Since the trace of $\mathbb{T}(r_1)$ is $K^2(r_1)$, the largest eigenvalue of $\mathbb{T}(r_1)$ is no smaller than $\frac{1}{\dim(\mathbb{V})} K^2(r_1)$. Use $u$ in what follows to denote a length one eigenvector of $\mathbb{T}(r_1)$ with the largest of $\mathbb{T}(r_1)$'s eigenvalues.

Because $K_u^2(r_1) \geq c_0^{-1} K^2(r_1)$, the value of $N_u(r_1)$ is no greater than $c_0^{-1} N(r_1)$. (That this is so follows from their respective definitions in (2.8) and (2.2).) Therefore, $N_u(r_1)$ is no greater tha $c_0 \varepsilon^{1/4}$. The function $N_u$ also obeys (2.5), and as a consequence, it obeys

$$N_u(r) \leq \varepsilon^{1/8} \text{ on the whole of } [\varepsilon^{1/(8n)} r_1, r_1]$$

(4.1)

if $\varepsilon \leq c_0^{-1}$. This implies in turn (because of (2.8)) that

$$K_u(r) \geq c_0^{-1} K(r_1) \text{ on the whole of } [\varepsilon^{1/(8n)} r_1, r_1].$$

(4.2)

Since $K_u(r) \leq K(r_1)$ in any event, the constant $K(r_1)$ can be used as a proxy for the function $K_u(\cdot)$ on the interval $[\varepsilon^{1/(8n)} r_1, r_1]$.

Let $v$ denote a length one eigenvector of $\mathbb{T}(r_1)$ with the smallest of $\mathbb{T}(r_1)$'s eigenvalues. The vectors $v$ and $u$ are orthogonal if $\mathbb{T}(r_1)$ is not a multiple of the identity. If it is, take $v$ to be orthogonal to $u$. (In fact, it can be shown that the minimum eigenvalue of $\mathbb{T}(r_1)$ is much less than the maximum eigenvalue.) Since $v$'s eigenvalue is $K_v^2(r_1)$, the latter is no greater than $K_u^2(r_1)$. Lemma 2.3 says that $K_v^2 \geq c_0^{-1} K_v^2(r_1)$ on the interval $[\varepsilon^{1/(8n)} r_1, r_1]$. Lemma 2.3 also says that $N_v \leq \varepsilon^{1/4}$ on this same interval.

*Part 2*: Because $u$ and $v$ are eigenvectors of $\mathbb{T}(r_1)$ and $v$ is orthogonal to $u$, the value of $\mathbb{T}(r_1)_{ab}u^a v^b$ is zero. Let $\mathcal{P}$ denote the function on $[\frac{1}{2}r_1, r_1]$ given by the rule $\mathcal{P}(\cdot) = \mathbb{T}(\cdot)_{ab} u^a v^b$. Since $\mathcal{P}(r_1) = 0$, the norm of $\mathcal{P}$ at any $s \in [\frac{1}{2}r_1, r_1]$ is bounded by

$$|\mathcal{P}(s)| \leq \int_s^{r_1} |\tfrac{d\mathcal{P}}{dr}| \, dr \ .$$

(4.3)

Meanwhile, the definition of $\mathbb{T}$ in (2.9) leads to a formula for the derivative of $\mathcal{P}$,

$$\frac{d\mathcal{P}}{dr} = \frac{1}{r^3} \int_{\partial B_r} (\langle a(u) \nabla_{Ar} a(v) \rangle + \langle \nabla_{Ar} a(u) a(v) \rangle) \ ,$$

(4.4)

which can be written using Stokes' theorem and (2.6) (and its $u$ analog) as

$$\frac{d\mathcal{P}}{dr} = \frac{2}{r^3} \int_{B_r} (\langle \nabla_A a(u) \nabla_A a(v) \rangle + \langle [a, a(u)], [a, a(v)] \rangle) \ .$$

(4.5)

The preceding identity leads to the bound $|\frac{d\mathcal{P}}{dr}| \leq \frac{2}{r} K_u(r) K_v(r) (N_u(r) N_v(r))^{1/2}$. And, given what was said in Part 1, this bound and (4.3) imply that

$$|\mathcal{P}(s)| \leq c_0 \, \varepsilon^{1/4} K_u(r_1) K_v(r_1) \ .$$

(4.6)

The take-away is that $|\mathcal{P}|$ is smaller than $K_u(r_1) K_v(r_1)$ on $[\frac{1}{2} r_1, r_1]$ by an $\mathcal{O}(\varepsilon^{1/4})$ factor.

*Part 3*: If $\varepsilon < c_0^{-1}$ (which will be assumed henceforth), then the interval $[\frac{1}{2}r_1, \frac{7}{8}r_1]$ is contained inside $[\varepsilon^{1/(8n)} r_1, r_1]$; the significance being that the latter interval appears in the assertions of Proposition 3.1. Because of the bound $N_u \leq \varepsilon^{1/8}$ on $[\frac{1}{2}r_1, \frac{7}{8}r_1]$, there is, by necessity, a subset of measure at least $\frac{1}{4}(1 - c_0 \varepsilon^{1/16}) r_1$ in the interval $[\frac{5}{8}r_1, \frac{7}{8}r_1]$ that is characterized as follows: If $s$ is in this subset, then

$$\int_{\partial B_s} (|\nabla_A a(u)|^2 + \|[a, a(u)]\|^2) \leq c_0 \varepsilon^{1/16} r_1^{n-3} K_u^2(r_1) \ .$$

(4.7)

Likewise, because $N_v(r_1) \leq \varepsilon^{1/8}$ on $[\frac{1}{2}r_1, \frac{7}{8}r_1]$, there is another a subset of measure at least $\frac{1}{4}(1 - c_0 \varepsilon^{1/16}) r_1$ in the interval $[\frac{5}{8}r_1, \frac{7}{8}r_1]$ with any given member $s$ obeying

$$\int_{\partial B_s} (|\nabla_A a(v)|^2 + \|[a, a(v)]\|^2) \leq c_0 \varepsilon^{1/16} r_1^{n-3} K_v^2(r_1) .$$

(4.8)

If $\varepsilon < c_0^{-1}$, then these two subsets of $[\frac{5}{8} r_1, \frac{7}{8} r_1]$ will intersect because their measures are nearly the full measure of the interval (which is $\frac{1}{4} r_1$). Therefore, supposing henceforth that $\varepsilon < c_0^{-1}$, choose a point $s \in [\frac{5}{8} r_1, \frac{7}{8} r_1]$ where both (4.7) and (4.8) hold.

*Part 4*: Let $\wp$ denote the function on $\partial B_s$ given by $\langle a(u) a(v) \rangle$. Its integral over $\partial B_s$ is $s^{n-1} \mathcal{P}(s))$. Since $|d\wp| \leq |a(u)||\nabla_A a(v)| + |\nabla_A a(u)||a(v)|$; and since $|a(u)| \leq c_0 K_u(r_1)$ (the top bullet of Lemma 2.3) and likewise $|a(v)| \leq c_0 K_v(r_1)$, it follows as a consequence of (4.7) and (4.8) that

$$\int_{\partial B_s} |d\wp|^2 \leq c_0 \varepsilon^{1/16} r_1^{n-3} K_u^2(r_1) K_v^2(r_1) .$$

(4.9)

Keeping in mind that the average value of $\wp$ is $\frac{1}{\omega} \mathcal{P}$, the bound in (4.9) (and the fact that the kernel of $d^\dagger d$ on $\partial B_s$ is spanned by the constant functions) leads to:

$$\int_{\partial B_s} |\wp - \tfrac{1}{\omega} \mathcal{P}|^2 \leq c_0 \varepsilon^{1/16} r_1^{n-1} K_u^2(r_1) K_v^2(r_1) .$$

(4.10)

And, what with (4.6), this leads in turn to an $L^2$ bound on $\langle a(u) a(v) \rangle$ (which is $\wp$):

$$\int_{\partial B_s} \langle a(u) a(v) \rangle^2 \leq c_0 \varepsilon^{1/16} r_1^{n-1} K_u^2(r_1) K_v^2(r_1) .$$

(4.11)

The take-away from all of this is that the integral of $\langle a(u)a(v) \rangle^2$ on $\partial B_s$ smaller than $r_1^{n-1} K_u^2(r_1) K_v^2(r_1)$ by an $\mathcal{O}(\varepsilon^{1/16})$ factor.

*Part 5*: The square of the norm of $[a, a(v)]$ is no smaller than $\|[a(u), a(v)]\|^2$ because $u$ and $v$ are orthogonal in $\mathfrak{su}(2)$. As a consequence, the contribution to the integrand in (4.8) from $\|[a, a(v)]\|^2$ leads to a bound on the integral over $\partial B_s$ of $\|[a(u), a(v)]\|^2$. Moreover, since

$$\|[a(u), a(v)]\|^2 = 4 |a(u)|^2 |a(v)|^2 - 4\langle a(u)a(v) \rangle^2 ,$$

(4.12)

(which is a property of the Lie algebra $\mathfrak{su}(2)$), the bound in (4.8) says, in part, that

$$\int_{\partial B_s} |a(u)|^2 |a(v)|^2 - \int_{\partial B_s} \langle a(u) a(v) \rangle^2 \leq c_0 \varepsilon^{1/16} r_1^{n-3} K_v^2(r_1) \,.$$

(4.13)

Together, the inequalities in (4.11) and (4.13) imply:

$$\int_{\partial B_s} |a(u)|^2 |a(v)|^2 \leq c_0 \varepsilon^{1/16} r_1^{n-1} K_u^2(r_1) K_v^2(r_1) \left(1 + \frac{1}{r_1^2 K_u^2(r_1)}\right).$$

(4.14)

Note with regards to the right hand side of (4.14) that $r_1^2 K_u^2(r_1) \geq r_1^2 K^2(1) \geq 1$ if $\varepsilon < c_0^{-1}$ because of the assumption that $\rho^\varepsilon K(1) > 1$ and the fact that $r_1 \geq \rho^{(1-30\sqrt{\varepsilon})}$. As a consequence, the right hand side of (4.14) is no greater than $c_0 \varepsilon^{1/16} r_1^{n-1} K_u^2(r_1) K_v^2(r_1)$.

The take-away from this part of the proof is that the integral of $|a(u)|^2 |a(v)|^2$ on $\partial B_s$ smaller than $r_1^{n-1} K_u^2(r_1) K_v^2(r_1)$ by an $\mathcal{O}(\varepsilon^{1/16})$ factor.

*Part 6*: As explained directly, the $r = r_1$ and vector $u$ version of Proposition 3.1 can be brought to bear with regards to the pointwise norm of $|a(u)|$ on the sphere $\partial B_s$. To this end, keep in mind that $N_u(s) \leq \varepsilon^{1/8}$ (which is implied by (4.1)); and that $K_u(s)$ is greater than $c_0^{-1} K_u(r_1)$ (which is implied by (4.2)), but that $K_u(s) \leq K_u(r_1)$ in any event (because $K_u$ is an increasing function). To exploit Proposition 3.1, introduce by way of notation $\Omega$ to denote the set of points in $\partial B_s$ where $|a(u)| \leq \frac{1}{2\sqrt{\omega}} K_u(r_1)$. The inequality in (4.14) holds with the integration domain restricted to the smaller domain $\partial B_s - \Omega$; and so (4.14) it leads to:

$$\int_{\partial B_s - \Omega} |a(v)|^2 \leq c_0 \varepsilon^{1/6} r_1^{n-1} K_v^2(r_1) \,.$$

(4.15)

Meanwhile, $|a(v)|$ on $\partial B_s$ is no greater than $c_0 K_v^2(r_1)$ (by virtue of the $r = r_1$ version of Proposition 3.1); and this bound on $\Omega$ implies that:

$$\int_\Omega |a(v)|^2 \leq c_0 r_1^{n-1} f_\Omega K_v^2(r_1)$$

(4.16)

with $f_\Omega$ denoting the fraction of the n-1 dimensional volume of $\partial B_s$ occupied by $\Omega$. The number $f_\Omega$ is no greater than $c_0 \sqrt{N_u(r_1)}$ (by virtue of the third bullet of the $r = r_1$ and vector $u$ version of Proposition 3.1) which is, therefore, at most $c_0 \varepsilon^{1/16}$ according to (4.1). Therefore, (4.15) and (4.16) together say that

$$\int_{\partial B_s} |a(v)|^2 \leq c_0 \varepsilon^{1/6} r_1^{n-1} K_v^2(r_1) \ ;$$

(4.17)

which is to say that $K_v^2(s) \leq c_0 \varepsilon^{1/16} K_v^2(r_1)$. The latter bound and what is said at the end of Part 1 to the effect that $K_v^2(s) \geq c_0^{-1} K_v^2(r_1)$ lead to the inequality

$$K_v^2(r_1) \leq c_0 \varepsilon^{1/6} K_v^2(r_1) \ .$$

(4.18)

This inequality is patently nonsense if $\varepsilon < c_0^{-1}$; and it is the desired nonsense that proves Theorem 1.1.

## 5. The Kapustin-Witten equations: Proofs of Theorems 1.2 and 1.3

This section gives the proof of Theorems 1.2 and 1.3. The proofs use $(A, a)$ to denote a solution to the relevant version of the Kapustin-Witten equations. Keep in mind that Theorem 1.1 can be invoked for $(A, a)$ since (1.1) is obeyed by $(A, a)$.

*Proof of Theorem 1.2*: By virtue of Theorem 1.1, the assertion of Theorem 1.2 holds except perhaps if $a \wedge a$ is identically zero. So, assume that such is the case. Theorem 1.2 then follows by proving that $\nabla_A a$ is identically zero. Note in this regard that if $a \wedge a = 0$, then $a$ (where it is not zero) can be written as $a = t\sigma$ with $t$ being a map to $\mathbb{R}^4$ (which should be viewed as a 1-form) and with $\sigma$ being a map to $\mathfrak{su}(2)$ with $|\sigma| = 1$. This depiction of $a$ is used below.

The arguement that $\nabla_A a = 0$ for the cases where $\tau = 0$ is given below in four parts. (The argument borrows from Ben Mare's Ph.D. thesis [M].) The cases where $\tau = 1$ is obtained from the $\tau = 0$ case by reversing the orientation of $\mathbb{R}^4$ and changing $a$ to $-a$.

*Part 1*: In the event that $\tau = 0$, the anti-self dual part of $d_A a$ is zero (the second bullet of (1.3)); and this equation with the equation $d_A *a = 0$ (the third bullet of (1.3)) can hold only if $\nabla_A \sigma = 0$ and both $(dt)^-$ and $d*t$ are identically zero. (You can see this by projecting the equations along $\sigma$ and orthogonal to $\sigma$ keeping in mind that each component of $\nabla_A \sigma$ is orthogonal in $\mathfrak{su}(2)$ to $\sigma$.) Therefore, $\sigma$ is A-covariantly constant where $a \neq 0$. However, because A is a smooth connection, $\sigma$ is defined on the whole of $\mathbb{R}^4$ and A-covariantly constant everywhere. This can be proved with the help of [T3] which can be brought to bear because the zero locus of $a$ is the zero locus of $|t|$ and $t$ is an example of what is said in [T3] to be a $\mathbb{Z}/2$ harmonic spinor. In particular, Theorem 1.3 in [T3] asserts that the zero locus of a $\mathbb{Z}/2$ harmonic spinor (such as $t$) has Hausdorff dimension 2, and it implies that the complement of the zero locus in any given ball in $\mathbb{R}^4$ is path connected. Given this input from [T3], then the $\mathfrak{su}(2)$ valued function $\sigma$ can be defined where $t$ is zero by its parallel transport (via A) along paths to where $t$ is not zero. Since A is smooth and $\sigma$ is already parallel where $t \neq 0$, two such paths to the same $t = 0$ point produce the same value of $\sigma$.

With $\sigma$ defined everywhere (and with it being A-covariantly constant), the 1-form $t$ is therefore also defined everywhere as $\langle \sigma\, a \rangle$. Meanwhile, any solution $s$ (and $s = t$ in particular) to the linear, constant coefficient equation $(ds)^- = 0$ and $d*s = 0$ on $\mathbb{R}^4$ has the following property: Either it is a constant 1-form, or its analog of $\kappa^2$, which is the function on $(0, \infty)$ given by

$$r \to \frac{1}{r^{n-1}} \int_{\partial B_r} |s|^2 \ ,$$

(5.1)

grows faster than a non-zero multiple of r. (This dichotomy can be proved by writing $s$ using radial coordinates on $\mathbb{R}^4$ and using $S^3$ spherical harmonics.)

***Proof of Theorem 1.3***: Suppose that $(A, a)$ obeys a version of (1.3) with $\tau$ neither 0 nor 1, and suppose that $a \wedge a = 0$. Then, the pair

$$\hat{A} = A - \frac{(2\tau - 1)}{\tau(\tau - 1)} a \quad \text{and} \quad \hat{a} = \frac{(1 - 2\tau + 2\tau^2)}{\tau(\tau - 1)} a$$

(5.2)

obeys the $\tau = \frac{1}{2}$ version of (1.3). This is to say that

- $F_{\hat{A}} = *d_{\hat{A}} \hat{a}$ .
- $d_{\hat{A}} * \hat{a} = 0$.

(5.3)

This understood, the discussion that follows talks only about the $\tau = \frac{1}{2}$ version of (1.3).

To see about solutions to the $\tau = \frac{1}{2}$ Kapustin-Witten equations, first write the Euclidean coordinates of $\mathbb{R}^4$ as $(x_1, x_2, x_3, x_4)$. Let $\mathbb{R}^3$ denote the span of the first three coordinates. A pair $(\mathcal{A}, \Phi)$ of connection on the product SU(2) bundle over the $(x_1, x_2, x_3)$ version of $\mathbb{R}^3$ and section of the associated $\mathfrak{su}(2)$ vector bundle over this $\mathbb{R}^3$ is said to obey the SU(2) monopole equations when

$$F_{\mathcal{A}} = *d_{\mathcal{A}} \Phi$$

(5.4)

with $*$ denoting in this equation the Euclidean Hodge star on $\mathbb{R}^3$. The only solutions to this equation with $\lim_{|x|\to\infty} |\Phi| = 0$ are equivalent (via the standard action of the gauge group $C^\infty(\mathbb{R}^3; SU(2))$) to the pair where $\Phi$ is zero and $\mathcal{A}$ is the product connection. On the other hand, there are a plethora of gauge inequivalent solutions with $\lim_{|x|\to\infty} |\Phi| = 1$. See e.g. [T4], [JT], [D]. Of course, there is the one gauge equivalent class with $\nabla_\mathcal{A}\Phi$ being identically zero. However, all of the others have non-trivial $\nabla_\mathcal{A}\Phi$. In fact, the integral over $\mathbb{R}^3$ of $|\nabla_\mathcal{A}\Phi|^2$ is $4\pi$ times a positive integer (call it k), and the space of gauge equivalence classes of solutions where this integral is $4\pi k$ with $k \geq 1$ is a smooth manifold of dimension $4k-1$.

Now, suppose that $(\mathcal{A}, \Phi)$ is a solution to (5.4) with $\lim_{|x|\to\infty} |\Phi| = 1$. Let A denote the pull-back of the connection $\mathcal{A}$ via the projection from $\mathbb{R}^4$ to $\mathbb{R}^3$. This is to say that A is the connection $\mathcal{A}$, but now viewed as a connection on the product principle SU(2) bundle over $\mathbb{R}^4$. Set $a$ to be the $\mathfrak{su}(2)$-valued 1-form $a = \Phi dx_4$. Then the pair $(A, a)$ obeys the $\tau = \frac{1}{2}$ version of (1.3); and the function $\kappa$ that is depicted in (1.2) has $\sqrt{2}\pi$ limit as $r \to \infty$. Note that $a \wedge a = 0$ (as required by Theorem 1.1). In addition, if the integral of $|\nabla_\mathcal{A}\Phi|^2$ over $\mathbb{R}^3$ is non-zero, then the integral of $|\nabla_A a|^2$ over the ball of radius r centered at the origin in $\mathbb{R}^4$ will grow like a non-zero multiple of r for large r.

To prove the last assertion of the theorem (about $L^2$ curvature when $a \wedge a = 0$ and $\frac{1}{r}\kappa(r)$ has limit zero as $r \to \infty$), suppose for argument's sake that $(A, a)$ obeys $F_A = *d_A a$ on $\mathbb{R}^4$ with $a \wedge a = 0$. In this event, the curvature 2-form of the connection $A + a$ is self-dual (and that of $A - a$ is anti-self dual). Let $\hat{A} = A + a$. Since $\hat{A}$ has square integrable, self-dual curvature, it is the pull-back of a connection on $S^4$ with self-dual curvature by the inverse to the stereographic projection map (this is proved in [U?]). As a consequence, $|F_{\hat{A}}| \leq c_0 \frac{1}{|x|^4}$ for $|x| \geq 1$. Keeping this in mind, and noting that $\nabla_{\hat{A}} a = \nabla_A a$, it follows that $F_{\hat{A}} = 2(d_{\hat{A}} a)^+$ and thus that $|F_{\hat{A}}|^2$ is equal to $2 * \langle F_{\hat{A}} \wedge d_{\hat{A}} a \rangle$. Therefore, if $r > 1$,

$$\int_{B_r} |F_{\hat{A}}|^2 = 2 \int_{B_r} \langle F_{\hat{A}} \wedge d_{\hat{A}} a \rangle = 2 \int_{\partial B_r} \langle F_{\hat{A}} \wedge a \rangle$$

(5.5)

with the right most inequality obtained using Stokes' theorem and the fact that $d_{\hat{A}} F_{\hat{A}} = 0$. By virtue of the bound $|F_{\hat{A}}| \leq c_0 \frac{1}{r^4}$ on $\partial B_r$, the right most integral in (5.5) can not be greater than $c_0 \frac{1}{r}\kappa(r)$. If the lim-inf as $r \to \infty$ of $\frac{1}{r}\kappa(r)$ is zero, then $F_{\hat{A}}$ is zero so $\hat{A}$ is flat (and thus gauge equivalent to the product connection.) Meanwhile, if $F_{\hat{A}}$ is zero, then so is $(d_A a)^+$ thus $(dt)^+$ is zero also. Keeping in mind that $d*t = 0$, it follows that $t$ obeys the first order, elliptic system $(dt)^+ = 0$, $d*t = 0$. As was the case with the $(ds)^- = 0$, $d*s = 0$ equations, this system also implies that the $s = t$ version of the function in (5.1) must grow faster than a non-zero multiple of r. Since this last event is precluded by the assumptions in the theorem, $t$ must vanish also. Therefore, $a = 0$ and $F_{\hat{A}} = 0$ so the solution is the trivial solution.